\documentclass[12pt,reqno]{amsart}
\usepackage[a4paper,margin=28mm]{geometry}
\usepackage{amsmath,amssymb,amsthm,mathtools,booktabs,array}
\usepackage[colorlinks=true,linkcolor=blue,citecolor=blue,urlcolor=blue]{hyperref}
\allowdisplaybreaks[2]
\theoremstyle{plain}
\newtheorem{theorem}{Theorem}[section]
\newtheorem{lemma}[theorem]{Lemma}
\newtheorem{proposition}[theorem]{Proposition}
\newtheorem{corollary}[theorem]{Corollary}
\theoremstyle{definition}

\numberwithin{equation}{section}
\newcommand{\kk}{\mathbb Q(q)}
\newcommand{\Res}{\operatorname{Res}}
\newcommand{\Rad}{\operatorname{Rad}}
\newcommand{\HC}{\operatorname{HC}}

\newcommand{\diag}{\operatorname{diag}}
\newcommand{\eps}{\epsilon}
\newcommand{\Span}{\operatorname{span}}
\newcommand{\tens}{\otimes}
\newcommand{\id}{\operatorname{id}}
\title[Spherical functions of type FII]{Spherical functions of the quantum symmetric pair of type FII}
\author{Shaolong Han}
\address{Department of Mathematics, Shantou University, Shantou 515063, P.R. China}
\email{slhan@stu.edu.cn}
\begin{document}
\begin{abstract}
We prove Schl\"osser's conjecture relating quantum FII spherical functions to Askey--Wilson polynomials for equal left and right coideal parameters and determine the corresponding radial operator.
An explicit transport formula gives the shifted restrictions for arbitrary nonzero parameters.
\end{abstract}
\maketitle

\section{Introduction}

Zonal spherical functions of quantum symmetric pairs give
representation-theoretic realizations of Macdonald and Koornwinder
polynomials. Letzter identified the distinguished Weyl-invariant
spherical family with Macdonald polynomials for standard pairs
with reduced restricted root systems \cite{Let04}.
The exceptional FII pair has a nonreduced restricted root system
of rank one. Its Askey--Wilson identification is conjectured
in Schl\"osser's finite-type radial-part framework \cite{Sch26}.
The standard FII case also lies outside the identification
results of \cite{MS26}.

We prove Schl\"osser's FII conjecture for the same standard
coideal on both sides.
For every $m\geq0$, the shifted toral restriction of the spherical
function of highest weight $m\omega_4$ is proportional to the
monic degree-$m$ Askey--Wilson polynomial with base $q^2$ and
parameters $(q^{11},-q^7,q^5,-q)$.
Highest-weight normalization fixes the scalars, and we determine
the corresponding radial operator (Theorem~\ref{thm:main}).

Letzter's torus-equivariance describes parameter transport for
quantum spherical functions \cite{Let03}.
In the FII family, rescaling weight components gives an explicit
formula valid for arbitrary nonzero left and right parameters.
For the normalized choices $\pm q^5$, their relative sign multiplies
all four Askey--Wilson parameters (Corollary~\ref{cor:transport}).

The existence of a spherical vector in the $26$-dimensional
fundamental $F_4$-module is known \cite{Wat23}; we verify its
coordinates in Proposition~\ref{prop:fundamental}.
Relevant symmetry properties are established in \cite{Mee26}.
Our finite calculation is the Cartan projection in the tensor
square. The relevant weight spaces have dimensions $12$ and $28$;
the black Levi invariants in the corresponding Cartan weight
spaces have dimensions one and two.
Lowering words certify Cartan membership, and Gram projections
give the degree-two spherical function.

Harish--Chandra surjectivity provides a lift, but the general
radial support bound allows half-step translations \cite{Sch26}.
A grading gives a lift whose radial part has integral translation
support. Weyl symmetry and the support bound then leave at most
three shifts.
Finite-support uniqueness determines the coefficient functions
from the constant and first two nonconstant spherical functions
and their Harish--Chandra eigenvalues.
Comparison with the Askey--Wilson operator, using the full
Harish--Chandra spectrum and distinct Askey--Wilson eigenvalues,
proves the identification in every degree.

The paper is organized as follows.
Section~\ref{sec:conventions} fixes conventions and states the conjecture and main result.
Section~\ref{sec:cartan} constructs spherical coefficients and
proves their shifted symmetry.
Sections~\ref{sec:smallspaces} and~\ref{sec:second} give Cartan
certificates and compute the degree-two function, respectively.
Section~\ref{sec:reconstruction} reconstructs the radial operator
and proves the main result.
Section~\ref{sec:transport} gives parameter transport, sign
dependence, and normalization at the identity
for the standard spherical functions.

\section{Conventions, conjecture and main result}
\label{sec:conventions}
We work over $\kk$, with $q$ indeterminate. Fix $I=\{1,2,3,4\}$ and $(e_i,e_j)=2\delta_{ij}$, $\alpha_i=e_{i+1}-e_{i+2}$ $(1\leq i\leq 2)$, $\alpha_3=e_4$, $\alpha_4=\tfrac12(e_1-e_2-e_3-e_4)$.
Let $\mathsf P$ be the weight lattice with fundamental weights
$\omega_i$, and $\mathsf Y$ its dual. For
$\nu\in\mathsf P\otimes\mathbb Q$, write $K_\nu=K_{h_\nu}$,
where $h_\nu\in\mathsf Y\otimes\mathbb Q$ and
$\langle h_\nu,\mu\rangle=(\nu,\mu)$ $(\mu\in\mathsf P)$,
using scalar and toral extensions when necessary.
On weight $\lambda\in\mathsf P$, $K_\nu$ acts by
$K_\nu(\lambda)=q^{(\nu,\lambda)}$.
Set $K_i=K_{\alpha_i}$, $q_i=q^{(\alpha_i,\alpha_i)/2}$, $q_1=q_2=q^2$, $q_3=q_4=q$ and $[n]=\frac{q^n-q^{-n}}{q-q^{-1}}$ $(n\in\mathbb Z)$.
Put
\[
\beta=\omega_4=e_1,\quad
\rho=\tfrac12(11e_1+5e_2+3e_3+e_4),\quad
(\beta,\beta)=2,\quad(\rho,\beta)=11.
\]

Let $U=U_q(F_4)$ be the quantum group of type $F_4$.
For Weyl group elements $w$, we use Lusztig's braid operators
$T_w=T''_{w,1}$ on $U$ and its integrable modules
\cite[37.1.3, 39.4.4 and~39.4.7]{Lus}.
Our Hopf conventions are
\begin{equation}\label{eq:coproduct}
		\Delta(K_h)=K_h\tens K_h,\quad
		\Delta(E_i)=E_i\tens1+K_i\tens E_i,\quad
		\Delta(F_i)=F_i\tens K_i^{-1}+1\tens F_i.
\end{equation}
The counit and antipode satisfy
\[
\begin{gathered}
	\eps(E_i)=\eps(F_i)=0,\  \eps(K_h)=1,\ 
	\mathsf S(E_i)=-K_i^{-1}E_i,\ 
	\mathsf S(F_i)=-F_iK_i,\ 
	\mathsf S(K_h)=K_{-h}.
\end{gathered}
\]
Contravariant forms are symmetric bilinear forms satisfying
\begin{equation}\label{eq:contravariance}
		(xu,v)=(u,\varrho(x)v),\quad
		\varrho(E_i)=q_i^{-1}F_iK_i,\quad
		\varrho(F_i)=q_iK_i^{-1}E_i,\quad
		\varrho(K_h)=K_h,
\end{equation}
where $\varrho$ is an algebra antiautomorphism.
The identity $\Delta\varrho=(\varrho\tens\varrho)\Delta$ makes
product forms contravariant.

For dominant $\lambda\in\mathsf P$, let $V(\lambda)$ be the simple
$U$-module of highest weight $\lambda$.
Finite-dimensional type-$1$ modules over $U$ and its Levi subalgebras
are semisimple, and their simple constituents have classical characters
\cite[Proposition~6.3.4(b) and Theorem~33.1.3(d)]{Lus}.

Let $s_i$ be the simple reflection for $\alpha_i$, and $w_X$ the
longest element of $\langle s_i:i\in X\rangle$, $X=\{1,2,3\}$.
Put $\Theta=-w_X$, acting dually by
$\langle\Theta h,\mu\rangle=\langle h,\Theta^{-1}\mu\rangle$
on $\mathsf Y$ and $\mathbb Q$-linearly on $\mathsf Y\otimes\mathbb Q$.
Set $\mathsf Y^\Theta=\{h\in\mathsf Y:\Theta h=h\}$ and
\[
U_X=\langle E_i,F_i,K_i^{\pm1}:i\in X\rangle,\quad
U^{0,\Theta}=\Span\{K_h:h\in\mathsf Y^\Theta\}.
\]
The FII coideal is $B=\langle U_X,U^{0,\Theta},B_4\rangle$, where
\begin{equation}\label{eq:B4}
		B_4=F_4+q^5E_\gamma K_4^{-1},\quad
		E_\gamma=T_{w_X}(E_4),\quad
		\gamma=w_X\alpha_4=\tfrac12(e_1+e_2+e_3+e_4).
\end{equation}
We follow \cite[Definition~2.4]{Sch26}.
Proposition~5.2 of \cite{Kol14}, with $c_4=-q^5$ and $s_4=0$
in Kolb's notation, gives the right coideal property. Adjoining fixed
group-like toral elements preserves it, so $\Delta(B)\subset B\tens U$.

For $\sigma\in\kk^\times$, put
\[
B_\sigma
=\langle U_X,U^{0,\Theta},F_4+\sigma E_\gamma K_4^{-1}\rangle,
\quad B=B_{q^5}.
\]
Let $\rho_X,\rho_X^\vee$ be the half sums of positive black roots
and coroots. The parameter conditions in
\cite[Definition~2.4]{Sch26} reduce to $\sigma^2=q^{10}$, since $(\alpha_4,\Theta\alpha_4-2\rho_X)=10$ and $\langle2\rho_X^\vee,\alpha_4\rangle=-6$.

Both choices $\sigma=\pm q^5$ satisfy
$\overline{\sigma}=\sigma^{-1}$ under the involution
$q\mapsto q^{-1}$ and may be chosen independently on the two
sides, as in the setup of \cite[Definition~2.8]{Sch26}.

For a finite-dimensional type-$1$ module $M$, write $M_\mu$ for
its $\mu$-weight space and
\[
M^A=\{v\in M:av=\eps(a)v\ (a\in A)\}
\]
for $A=B_\sigma,U_X$.
For $\delta,\sigma\in\kk^\times$, a
$(B_\delta,B_\sigma)$-spherical coefficient has right vector
$v\in M^{B_\sigma}$ and a left covector $\ell\in M^*$ satisfying
$\ell(bu)=\eps(b)\ell(u)$ ($b\in B_\delta$, $u\in M$).
The coefficient $c_{\ell,v}(x)=\ell(xv)$ ($x\in U$) is elementary
when it lies in one simple module's coefficient space.

The restricted root system and spherical lattice are
\begin{equation}\label{eq:sphericallattice}
		\Sigma=\{\pm\beta/2,\pm\beta\},\ 
		2L=\left\{\mu\in\mathsf P:\Theta\mu=-\mu,\ 
		\frac{(\alpha,\mu)}{(\alpha,\alpha)}\in\mathbb Z
		\text{ for all }\alpha\in\Sigma\right\}
		=\mathbb Z\beta.
\end{equation}
$\Theta\mu=-\mu$ forces $\mu=s\beta$, with integrality
conditions $s,2s\in\mathbb Z$.
Define $2\widetilde L$ by the same conditions as $2L$ in
\eqref{eq:sphericallattice}, with $\mathsf P$ replaced by
$\mathsf P\otimes\mathbb Q$.
The same integrality conditions give $2\widetilde L=2L$.

Write $e^\mu$ for the standard basis elements of $\kk[\mathsf P]$.
For a weight decomposition $v=\sum_\mu v_\mu$, set
\begin{equation}\label{eq:restriction}
	\Res'c_{\ell,v}=\sum_\mu\ell(v_\mu)e^\mu,\quad
	\Res c_{\ell,v}=\sum_\mu q^{-(\rho,\mu)}\ell(v_\mu)e^\mu.
\end{equation}
The second restriction uses the $-\rho$ shift of
\cite[Definition~2.12]{Sch26}.
Put $z=e^\beta$ and $S=z+z^{-1}$.
For $m\in\mathbb Z_{\ge0}$ and $\delta,\sigma\in\kk^\times$,
write $\phi_m^{\delta,\sigma}$ for the elementary
$(B_\delta,B_\sigma)$-spherical coefficient of $V(m\beta)$ with
right highest-weight coefficient and left covector value both $1$,
relative to the chosen highest vector. Put
$\phi_m=\phi_m^{q^5,q^5}$.
Uniqueness for $\delta,\sigma\in\{q^5,-q^5\}$ follows from
\cite[Lemma~4.2]{Mee26}.

For the structural results, extend scalars to
$\mathbb F=\overline{\mathbb C(q)}$ and adjoin all $K_h$,
$h\in\mathsf Y\otimes\mathbb Q$, with the usual weight action
and compatible rational powers $q^r$.
Denote this algebraic torus extension by $\widetilde U$;
its elements are finite PBW sums. Over $\mathbb F$, set
\[
\widetilde B
=\langle B,K_h:h\in(\mathsf Y\otimes\mathbb Q)^\Theta\rangle.
\]
Using $\Delta(x)=\sum x_{(1)}\tens x_{(2)}$, put
\[
\operatorname{ad}_{\rm r}(x)y
=\sum\mathsf S(x_{(1)})yx_{(2)},\quad
\mathcal N^-=\langle\operatorname{ad}_{\rm r}(U_X)(F_4)\rangle.
\]
Under the $\mathbb F$-linear isomorphism induced by multiplication
\cite[proof of Lemma~4.8]{Sch26}, the Harish--Chandra map is
\begin{equation}\label{eq:HCdefinition}
	\widetilde U\simeq
	\mathcal N^-\tens\mathbb F[K_{\mathbb Q\beta}]
	\tens\widetilde B,\quad
	\HC=\eps\tens\id\tens\eps.
\end{equation}
Here $\mathbb F[K_{\mathbb Q\beta}]$ is the group algebra generated
by $K_{s\beta}$, $s\in\mathbb Q$, and $\HC$ is the unshifted
Harish--Chandra map of \cite[Definition~3.4]{Sch26}.

Type-$1$ modules and matrix coefficients extend to
$\widetilde U$. Write $Z_{\widetilde U}(\widetilde B)$ for the
centralizer of $\widetilde B$ in $\widetilde U$.
We write $\Rad$ for the shifted radial part compatible with $\Res$.
For a $(B,B)$-spherical coefficient $\phi$ and
$x\in Z_{\widetilde U}(\widetilde B)$,
\begin{equation}\label{eq:radialaction}
	\Res(x\triangleright\phi)=\Rad(x)\Res\phi,\quad
	(x\triangleright\phi)(y)=\phi(yx),
\end{equation}
by \cite[Definition~3.2 and Lemma~3.3]{Sch26}.
The coefficient field is initially
$\mathcal F=\operatorname{Frac}(\mathbb F[\mathsf P])\supset\kk(z)$.
For $h\in\mathsf Y\otimes\mathbb Q$, evaluate $\mathcal F$ by
$e^\mu(h)=q^{\langle h,\mu\rangle}$ wherever denominators are nonzero.
For $\nu\in\mathsf P\otimes\mathbb Q$, define the
$\mathbb F$-automorphism $T_\nu$ of $\mathcal F$ by
\[
T_\nu e^\mu=q^{(\nu,\mu)}e^\mu,\quad
T_\beta g(z)=g(q^2z).
\]

Fix
\begin{equation}\label{eq:AWparameters}
	(a,b,c,d)=(q^{11},-q^7,q^5,-q).
\end{equation}
The Askey--Wilson operator and its eigenvalues are
\begin{align}
	D_{\rm AW}
	&=A_0(z)(T_\beta-1)+A_0(z^{-1})(T_{-\beta}-1),
	\label{eq:DAW}\\
	A_0(z)
	&=\frac{(1-az)(1-bz)(1-cz)(1-dz)}
	{(1-z^2)(1-q^2z^2)},
	\label{eq:AWcoefficient}\\
	\delta_m
	&=(q^{-2m}-1)(1-abcd\,q^{2m-2})
	\qquad(m\in\mathbb Z_{\ge0}).
	\label{eq:AWeigenvalues}
\end{align}

For $g\in\kk[S]$, \eqref{eq:DAW} shows that $D_{\rm AW}g$
is symmetric under $z\mapsto z^{-1}$.
The possible simple poles at $z=\pm1$ cancel between the two terms.
At $z^2=q^{-2}$, $g(q^2z)=g(z^{-1})=g(z)$; at $z^2=q^2$,
$g(q^{-2}z)=g(z^{-1})=g(z)$.
Thus $D_{\rm AW}g$ has no finite nonzero poles and is a
symmetric Laurent polynomial, hence belongs to $\kk[S]$.
Expansion of \eqref{eq:DAW} at $z=\infty$ gives
\begin{equation}\label{eq:AWtriangular}
	D_{\rm AW}S^m-\delta_mS^m
	\in\Span\{S^j:0\le j<m\}
	\quad(m\ge0).
\end{equation}

Equations~\eqref{eq:AWparameters} and~\eqref{eq:AWeigenvalues}
give
\begin{equation}\label{eq:deltam-deltan}
	\delta_m-\delta_n
	=q^{22}(q^{2m}-q^{2n})(1-q^{-22-2m-2n})\ne0
	\quad(m>n\ge0).
\end{equation}
Equations~\eqref{eq:AWtriangular} and~\eqref{eq:deltam-deltan}
therefore determine, for each $m\ge0$, a unique monic polynomial
$P_m(z)\in\kk[S]$ of degree $m$ satisfying
$D_{\rm AW}P_m=\delta_mP_m$.
These are the monic symmetric Laurent Askey--Wilson polynomials,
with $P_0=1$. For $m=1$,
\eqref{eq:DAW}--\eqref{eq:AWeigenvalues} give
\begin{equation}\label{eq:firstAW}
	P_1=S+\frac{a+b+c+d-abc-abd-acd-bcd}{abcd-1}
	=S+\frac{[2]}{[3]}.
\end{equation}
Negating all four parameters in \eqref{eq:AWparameters}
preserves $abcd$ and hence the eigenvalues in
\eqref{eq:AWeigenvalues}.
The same construction defines
$P_m(z;\varepsilon a,\varepsilon b,\varepsilon c,\varepsilon d\mid q^2)$
for $\varepsilon\in\{1,-1\}$.

With the same standard coideal $B$ on both sides,
Schl\"osser's FII conjecture
\cite[Lemma~6.3 and Conjecture~6.4]{Sch26} specializes to
\begin{equation}\label{eq:Schconjecture}
 \Res\phi_m\in\kk^\times P_m(z)
 \ (m\in\mathbb Z_{\ge0}).
\end{equation}
Here the source's $q_\beta$ equals $q$, and its factor
$e^{\mu-m\beta}$ is $1$ because $\mu=m\beta$.

Let $R_s=\{\pm e_i:1\leq i\leq4\}
\cup\{\tfrac12(\pm e_1\pm e_2\pm e_3\pm e_4)\}$.
The quasi-minuscule model of $V=V(\beta)$ has basis
$x_\mu$ ($\mu\in R_s$), $y_3,y_4$, of weights $\mu,0,0$.
Its neighboring-root actions are
\[
E_ix_\mu=x_{\mu+\alpha_i}\ (\mu+\alpha_i\in R_s),\quad
F_ix_\mu=x_{\mu-\alpha_i}\ (\mu-\alpha_i\in R_s).
\]
The remaining nonzero root-operator actions are
\begin{align}
	E_ix_{-\alpha_i}&=y_i,&
	F_ix_{\alpha_i}&=y_i,&
	E_iy_i&=[2]x_{\alpha_i},&
	F_iy_i&=[2]x_{-\alpha_i}
	&&(i=3,4),\label{eq:zerostrings}\\
	E_3y_4&=x_{\alpha_3},&
	F_3y_4&=x_{-\alpha_3},&
	E_4y_3&=x_{\alpha_4},&
	F_4y_3&=x_{-\alpha_4}.
	\label{eq:crosszero}
\end{align}
All unspecified entries vanish, and $K_ix_\mu=q^{(\alpha_i,\mu)}x_\mu$, $K_iy_j=y_j$.

In the following rank-two restrictions, bold numerals denote
irreducible modules of the indicated dimensions, and the
superscript $*$ denotes the dual.

On the color pairs $(1,2),(2,3),(3,4)$, the given operators yield
\[
3\mathbf3\oplus3\mathbf3^*\oplus8\mathbf1,\qquad
\mathbf5\oplus4\mathbf4\oplus5\mathbf1,\qquad
\mathbf8\oplus3\mathbf3\oplus3\mathbf3^*,
\]
respectively, with standard $A_2,B_2,A_2$ actions; the adjoint
$\mathbf8$ realizes \eqref{eq:crosszero}.
These rank-two actions and commutativity of nonadjacent colors
give a finite-dimensional type-$1$ $U$-module.
Since highest weights are dominant,
\[
\begin{gathered}
	\{\mu\in R_s:(\alpha_i,\mu)\ge0\text{ for all }i\in I\}
	=\{\beta\},\\
	\det\begin{pmatrix}[2]&1\\1&[2]\end{pmatrix}=[3]\ne0
\end{gathered}
\quad\Longrightarrow\quad
\bigcap_{i\in I}\ker E_i=\kk x_\beta.
\]
Complete reducibility identifies the model with $V(\beta)$.

Normalize the contravariant form by $(x_\beta,x_\beta)=1$.
Weight spaces are orthogonal, and
\begin{equation}\label{eq:GramV}
	(x_\mu,x_\nu)=\delta_{\mu\nu},\quad
	\bigl((y_i,y_j)\bigr)_{i,j=3,4}
	=q^{-1}\begin{pmatrix}[2]&1\\1&[2]\end{pmatrix}.
\end{equation}
Set
\begin{equation}\label{eq:h}
	h=y_4-[2]^{-1}y_3,\quad
	(h,h)=q^{-1}\frac{[3]}{[2]}.
\end{equation}
Then $(h,y_3)=0$ and
$V^{U_X}=\Span\{x_\beta,h,x_{-\beta}\}$.

Set
\begin{equation}\label{eq:w}
 w=x_\beta-q^{-5}\frac{[2]}{[3]}h+q^{-11}x_{-\beta}.
\end{equation}

\begin{theorem}\label{thm:main}
For every $m\in\mathbb Z_{\ge0}$,
\begin{equation}\label{eq:main}
 \Res\phi_m=q^{-11m}P_m(z).
\end{equation}
There exists $x\in Z_{\widetilde U}(\widetilde B)$ with
$\HC(x)=K_\beta+q^{-22}K_{-\beta}$ and
\begin{equation}\label{eq:mainradial}
 \Rad(x)=q^{-22}D_{\rm AW}+(1+q^{-22})\id.
\end{equation}
\end{theorem}

Theorem~\ref{thm:main} will be proved at the end of
Section~\ref{sec:reconstruction}.

\section{Spherical coefficients and Cartan projections}
\label{sec:cartan}

Define
\begin{equation}\label{eq:f}
 f=(w,\,\cdot\,),\quad
 f(x_\beta)=1,\quad f(y_4)=-q^{-6},\quad f(y_3)=0,\quad
 f(x_{-\beta})=q^{-11}.
\end{equation}
The values follow from \eqref{eq:GramV}; $f$ vanishes on the
other root vectors. Put $\phi_1=c_{f,w}$.

\begin{proposition}\label{prop:fundamental}
	The vector $w$ and covector $f$ are spherical, with highest-weight
	coefficient and value $1$, respectively. Moreover,
	\begin{equation}\label{eq:firstsphere}
				q^{11}\Res\phi_1=z+z^{-1}+\frac{[2]}{[3]}.
	\end{equation}
\end{proposition}

\begin{proof}
The reduced expression
$w_X=(s_1s_2s_3s_2s_1)(s_2s_3s_2)s_3$
gives effective spin-string colors $3,2,3,1,2,3$.
The upward braid coefficients are $1$, and the reverse product is
$(-1)^6q_3^3q_2^2q_1=q^9$. Hence \eqref{eq:B4} gives
\[
\begin{aligned}
 B_4x_\beta=x_\gamma,\quad
 B_4x_{-\beta}=q^6x_{-\alpha_4},\quad
 B_4h=\frac{[3]}{[2]}(x_{-\alpha_4}+q^5x_\gamma).
\end{aligned}
\]
Together with \eqref{eq:w}, these formulas give $B_4w=0$.
Since $w\in V^{U_X}$ and $\langle\mathsf Y^\Theta,\beta\rangle=0$,
the right spherical conditions hold.
The same braid computation gives
\begin{equation}\label{eq:whiteimages}
 B_4x_{\alpha_4}=y_4+q^{-6}x_\beta,\quad
 B_4x_{-\gamma}\equiv x_{-\beta}+q^{-5}y_4
 \pmod{\kk y_3}.
\end{equation}
The black and fixed-torus left conditions for $f$ follow from
\eqref{eq:contravariance}.
Weight support leaves only $x_{\alpha_4},x_{-\gamma}$ for $fB_4=0$.
Since $f(y_3)=0$, equations~\eqref{eq:whiteimages} and~\eqref{eq:f} give
\[
 f(B_4x_{\alpha_4})=-q^{-6}+q^{-6}=0,\quad
 f(B_4x_{-\gamma})=q^{-11}-q^{-5}q^{-6}=0.
\]
The normalization follows from \eqref{eq:w} and~\eqref{eq:f};
\eqref{eq:restriction} gives \eqref{eq:firstsphere}.
\end{proof}

For $m\in\mathbb Z_{>0}$, semisimplicity and the multiplicity-one
weight $m\beta$ identify $V(m\beta)$ with
$V_m=Ux_\beta^{\tens m}\subset V^{\tens m}$ by matching highest vectors.
The restricted contravariant form has $U$-stable radical.
Since $V_m$ is simple and
$(x_\beta^{\tens m},x_\beta^{\tens m})=1$, this form is nondegenerate.
Thus $V^{\tens m}=V_m\oplus V_m^\perp$; both summands are
$U$-stable by \eqref{eq:contravariance}, so the orthogonal projection
$\Pi_m$ is $U$-linear.

Put $V_0=\kk$, $\Pi_0=\id$, $w^{\tens0}=f^{\tens0}=1$, and
\begin{equation}\label{eq:allspheres}
	w_m=\Pi_m(w^{\tens m}),\quad
	f_m=(f^{\tens m})|_{V_m}=(w_m,\cdot),\quad
	\phi_m=c_{f_m,w_m}.
\end{equation}
In particular, $\phi_0=1$.

\begin{lemma}\label{lem:allspheres}
	For every $m\geq0$, $w_m$ and $f_m$ are spherical, with
	highest-weight coefficient and value $1$, respectively. Moreover,
	\[
	\Res'\phi_m-z^m
	\in\Span\{z^k:k\in\mathbb Z,\ -m\leq k<m\}.
	\]
\end{lemma}

\begin{proof}
For spherical vectors $u,v$, spherical covectors $\ell_1,\ell_2$,
and $b\in B$, $\Delta(b)\in B\tens U$ and
$(\eps\tens\id)\Delta(b)=b$ give
\[
\begin{aligned}
 \Delta(b)(u\tens v)=u\tens bv=\eps(b)(u\tens v),\quad
 (\ell_1\tens\ell_2)\Delta(b)
 =\ell_1\tens(\ell_2\circ b)
 =\eps(b)(\ell_1\tens\ell_2).
\end{aligned}
\]
Induction, Proposition~\ref{prop:fundamental}, and the $U$-linearity
	of $\Pi_m$ prove the spherical conditions.
	Since $\Pi_m$ preserves weights and fixes $x_\beta^{\tens m}$,
	\eqref{eq:w} gives the asserted support and highest-weight
	coefficient, while
	$f_m(x_\beta^{\tens m})=f(x_\beta)^m=1$.
	The restriction formula \eqref{eq:restriction} completes the proof.
\end{proof}

For dominant $\lambda$, \cite[Proposition~4.1(i)]{Mee26} gives
$0\ne v\in V(\lambda)^B\Rightarrow v_\lambda\ne0$.
The black torus conditions therefore imply
\[
K_iv=v\ (i\in X)
\quad\Longrightarrow\quad
(\alpha_i,\lambda)=0\ (i\in X)
\quad\Longrightarrow\quad
\lambda\in\mathbb Z_{\geq0}\beta.
\]
Together with Lemma~\ref{lem:allspheres}, this exhausts the spherical
highest weights.

Set $\mathsf{F}_m=q^{11m}\Res\phi_m$ for $m\geq0$.

\begin{lemma}\label{lem:shiftedsymmetry}
	For every $m\geq0$,
	$\Res\phi_m(z)=\Res\phi_m(z^{-1})$, and $\mathsf{F}_m$ is monic
	of degree $m$ in $S$.
\end{lemma}

\begin{proof}
	With $\operatorname{Ad}(a)(x)=axa^{-1}$ for invertible $a$,
	$\mathsf S^{-2}=\operatorname{Ad}(K_{2\rho})$ and
	\cite[Proposition~4.25]{Mee26} show that $\phi_m\circ\mathsf S$ is
	$(\operatorname{Ad}(K_{2\rho})B,B)$-spherical.
	So is $\psi(x)=\phi_m(K_{-2\rho}x)$, since $K_{-2\rho}(K_{2\rho}bK_{-2\rho})=bK_{-2\rho}$ $(b\in B)$.
	The longest Weyl element of $F_4$ is $-1$, so both functions
	belong to the same elementary block. The map
	$\ell\mapsto\ell\circ K_{-2\rho}$ identifies the left invariant
	covector spaces for $B$ and $\operatorname{Ad}(K_{2\rho})B$.
	Hence \cite[Lemma~4.2]{Mee26} makes the two functions proportional.
	
	Their common value at $K_\rho$ is nonzero.
	Indeed, at a generic positive real specialization of $q$,
	the Gram eigenvalues $1$ and $q^{-1}([2]\pm1)$ in
	\eqref{eq:GramV} are positive. Hence the product form is
	positive definite, and the positive eigenvalues of $K_{-\rho}$ give
	\[
	(\phi_m\circ\mathsf S)(K_\rho)
	=\psi(K_\rho)
	=(w_m,K_{-\rho}w_m)>0.
	\]
	Thus $\phi_m\circ\mathsf S=\psi$. Evaluating at $K_{\nu+\rho}$ gives $\phi_m(K_{-\nu-\rho})=\phi_m(K_{\nu-\rho})$ $(\nu\in\mathsf P\otimes\mathbb Q)$.
	Equation~\eqref{eq:restriction} yields the symmetry;
	Lemma~\ref{lem:allspheres} yields monicity and degree.
\end{proof}

\section{Small invariant spaces and finite certificates}
\label{sec:smallspaces}
Set $M=V\tens V$.
Set $A=q^2+1$, $B_0=q^4-q^2+1$,
$C=q^8-q^6+q^4-q^2+1$ and $D=\sum_{j=0}^6q^{2j}$.
For $\varepsilon,\delta\in\{\pm1\}$ and
$s=(s_2,s_3,s_4)\in\{\pm1\}^3$, put $n_j=(1-s_j)/2$ $(j=2,3,4)$ and
\begin{equation}\label{eq:OmegaSpin}
	\Omega_{\varepsilon\delta}
	=
	\sum_{s\in\{\pm1\}^3}
	(-1)^{n_2+n_4}q^{5n_2+3n_3+n_4}
	x_{(\varepsilon e_1+s_2e_2+s_3e_3+s_4e_4)/2}
	\tens
	x_{(\delta e_1-s_2e_2-s_3e_3-s_4e_4)/2}.
\end{equation}
Also set
\begin{align}\label{eq:Omega7}
	\Omega_7={}&
	x_{e_2}\tens x_{-e_2}
	-q^2x_{e_3}\tens x_{-e_3}
	+q^4x_{e_4}\tens x_{-e_4}
	-\frac{q^6}{[2]}y_3\tens y_3\notag\\
	&+q^6x_{-e_4}\tens x_{e_4}
	-q^8x_{-e_3}\tens x_{e_3}
	+q^{10}x_{-e_2}\tens x_{e_2}.
\end{align}
The vectors $\Omega_{\varepsilon\delta}$ and $\Omega_7$
are $U_X$-invariant.
In \eqref{eq:OmegaSpin}, raising the first spin weight by
$\alpha_i$ changes its coefficient by $-q_i^{-1}$;
the two contributions in \eqref{eq:coproduct} therefore cancel.
The reverse ratios give lowering invariance.
The same verification for \eqref{eq:Omega7} uses
\eqref{eq:zerostrings} at weight zero.

Write $M_\lambda^{U_X}=(M^{U_X})_\lambda$.
The black Levi restriction is
\[
V|_{U_X}
=
\mathbf7_0\oplus\mathbf8_{1/2}\oplus\mathbf8_{-1/2}
\oplus\mathbf1_1\oplus\mathbf1_0\oplus\mathbf1_{-1},
\]
where $\mathbf7,\mathbf8$ are the self-dual vector and spin
$B_3$-modules, and subscripts record the $e_1$ coefficient.
Schur's lemma gives $\dim M_\beta^{U_X}=3$ and
$\dim M_0^{U_X}=6$.

Ordered orthogonal bases of $M_\beta^{U_X}$ and $M_0^{U_X}$,
with Gram matrices $G_\beta$ and $G_0$, respectively, are:
\begin{equation}\label{eq:betabasis}
	p=x_\beta\tens h,\quad
	p'=h\tens x_\beta,\quad
	\Omega_{++},
	\quad G_\beta=\diag((h,h),(h,h),A^3B_0C),
\end{equation}
\begin{equation}\label{eq:zerobasis}
	p_+=x_\beta\tens x_{-\beta},\quad
	p_-=x_{-\beta}\tens x_\beta,\quad
	h\tens h,\quad
	\Omega_7,\quad\Omega_{+-},\quad\Omega_{-+},
\end{equation}
with
\begin{equation}\label{eq:Gzero}
	G_0=\diag(1,1,(h,h)^2,CD,A^3B_0C,A^3B_0C).
\end{equation}
Orthogonality and the Gram entries follow from
\eqref{eq:GramV}--\eqref{eq:h}.

At $q=1$, retain notation for classical modules.
Tensor interchange sends
$\Omega_{\varepsilon\delta}$ to $\Omega_{\delta\varepsilon}$
and fixes $\Omega_7$.
Together with the Weyl character decomposition
\[
\operatorname{Sym}^2V
=V(2\beta)\oplus V(\beta)\oplus V(0),
\]
this gives the invariant multiplicities
\begin{equation}\label{eq:invariantcounts}
	\dim V(2\beta)_\beta^{U_X}=1,\quad
	\dim V(2\beta)_0^{U_X}=2.
\end{equation}
The same dimensions hold for indeterminate $q$ by equality
of Weyl characters.

For $\lambda\in\{\beta,0\}$, let $\mathsf Q_\lambda$ be the
orthogonal projection of $M_\lambda$ onto $M_\lambda^{U_X}$.
Since $\Pi_2$ is a weight-preserving self-adjoint $U_X$-map,
it preserves $M_\lambda^{U_X}$ and its orthogonal complement
in $M_\lambda$. Hence
\begin{equation}\label{eq:commutingprojections}
	\Pi_2\mathsf Q_\lambda=\mathsf Q_\lambda \Pi_2.
\end{equation}

For a word $J=(i_1,\ldots,i_r)$, put
$F[J]=\Delta(F_{i_r})\cdots\Delta(F_{i_1})$;
letters act from left to right.
In the ordered bases \eqref{eq:betabasis}
and \eqref{eq:zerobasis}, define
\begin{equation}\label{eq:rbeta}
	r_\beta
	=
	\left(1,q^{-2},\frac{1}{q^2A^2B_0}\right)^{\mathsf T},
\end{equation}
\begin{equation}\label{eq:Rzero}
	R_0=
	\begin{pmatrix}
		1&q^{-4}&0&
		\dfrac{AB_0^2}{q^2CD}&
		-\dfrac{q^2-1}{qAC}&
		\dfrac{q^2-1}{q^4AC}\\[2mm]
		0&0&1&
		-\dfrac{q^5}{ACD}&
		\dfrac{q^{10}}{A^3B_0C}&
		\dfrac{1}{q^3A^3B_0C}
	\end{pmatrix}^{\mathsf T}.
\end{equation}

For $\lambda\in\{\beta,0\}$, let $Z_1,\ldots,Z_d$ denote the
ordered orthogonal basis of $M_\lambda^{U_X}$ in
\eqref{eq:betabasis} or \eqref{eq:zerobasis}, respectively.
For $v\in M_\lambda$, the $j$th coordinate of $\mathsf Q_\lambda v$ is
\begin{equation}\label{eq:Qcoordinate}
	\frac{(Z_j,v)}{(Z_j,Z_j)}
	\ (1\leq j\leq d),
\end{equation}
with pairings determined by \eqref{eq:GramV}.
Set
\begin{align*}
	I_1&=(4,4,3,2,1,3,2,3,4,3,2,1,3,2,3,4),\\
	I_2&=(4,3,2,3,4,1,2,3,4,3,2,1,3,2,3,4).
\end{align*}

\begin{proposition}\label{prop:smallspaces}
	The vector $r_\beta$ and the columns of $R_0$ form bases of
	$V(2\beta)_\beta^{U_X}$ and $V(2\beta)_0^{U_X}$, respectively.
\end{proposition}

\begin{proof}
	By \eqref{eq:commutingprojections}, $\mathsf Q_\lambda$ preserves
	$V(2\beta)_\lambda$.
	Equations~\eqref{eq:coproduct} and \eqref{eq:Qcoordinate},
	with the explicit action on $V$, give
	\begin{align}
		\mathsf Q_\beta F[4,3,2,1,3,2,3,4]x_\beta^{\tens2}
		&=r_\beta,\label{eq:wordbeta}\\
		\mathsf Q_0\bigl(F[I_1]x_\beta^{\tens2},F[I_2]x_\beta^{\tens2}\bigr)
		&=R_0
		\begin{pmatrix}
			0&1\\ q^{-3}A^2&2q^{-1}
		\end{pmatrix}.
		\label{eq:wordzero}
	\end{align}
	The coefficient matrix in \eqref{eq:wordzero} has determinant
	$-q^{-3}A^2\ne0$, so both columns of $R_0$ belong to
	$V(2\beta)_0^{U_X}$.
	By \eqref{eq:wordbeta}, $0\ne r_\beta\in V(2\beta)_\beta^{U_X}$.
	The first and third rows of $R_0$ give rank $2$, so
	\eqref{eq:invariantcounts} proves both assertions.
\end{proof}

\section{The second spherical function}\label{sec:second}

Once its translation support is contained in $\{-\beta,0,\beta\}$,
the radial operator will be determined by the shifted restrictions
of $\phi_0,\phi_1,\phi_2$ and their Harish--Chandra eigenvalues.
Set
\begin{equation}\label{eq:smallGram}
	H=R_0^{\mathsf T}G_0R_0,\quad
	\boldsymbol{\ell}=\binom{q^{-11}(1+q^{-4})}{q^{-12}}.
\end{equation}
By Lemma~\ref{lem:shiftedsymmetry}, we have
\begin{equation}\label{eq:secondsphere}
	\mathsf{F}_2=S^2+uS+v,\quad u,v\in\kk.
\end{equation}

\begin{proposition}\label{prop:second}
	The matrix $H$ is nonsingular, and the coefficients in
	\eqref{eq:secondsphere} satisfy
	\begin{equation}\label{eq:uclosed}
		u=\frac{qA^3B_0}{D},\quad
		v=\frac{[2]\bigl([2]-(q-q^{-1})^2[5][12]\bigr)}{[7][13]}.
	\end{equation}
\end{proposition}

\begin{proof}
	By \eqref{eq:allspheres},
	$\phi_2(x)=(f\tens f)\Delta(x)\Pi_2(w\tens w)$ for $x\in U$.
	Equation~\eqref{eq:w} gives
	\begin{equation}\label{eq:tensorweights}
			(w\tens w)_\beta
			=-q^{-5}\frac{[2]}{[3]}(p+p'),\quad
			(w\tens w)_0
			=q^{-10}\frac{[2]^2}{[3]^2}(h\tens h)
			+q^{-11}(p_++p_-).
	\end{equation}
	Use coordinate columns in the ordered invariant bases
	\eqref{eq:betabasis}--\eqref{eq:zerobasis}.
	At weight $\beta$, \eqref{eq:f}, \eqref{eq:h},
	\eqref{eq:rbeta}, and~\eqref{eq:tensorweights} give
	\begin{equation}\label{eq:betapairings}
			r_\beta^{\mathsf T}G_\beta r_\beta
			=\frac{D}{q^6AB_0},\quad
			(f\tens f)r_\beta
			=r_\beta^{\mathsf T}G_\beta(w\tens w)_\beta
			=-q^{-6}(1+q^{-2}).
	\end{equation}
	Proposition~\ref{prop:smallspaces}, \eqref{eq:restriction},
	and~\eqref{eq:betapairings} therefore give
	\[
	u=q^{11}(f\tens f)\Pi_2((w\tens w)_\beta)
	=\frac{q^{-1}(1+q^{-2})^2}
	{r_\beta^{\mathsf T}G_\beta r_\beta}
	=\frac{qA^3B_0}{D}.
	\]
	
	At weight zero, \eqref{eq:Rzero}, \eqref{eq:Gzero},
	\eqref{eq:smallGram}, and~\eqref{eq:tensorweights} give
	\begin{equation}\label{eq:zeropairings}
		(f\tens f)R_0=\boldsymbol{\ell}^{\mathsf T},\quad
		R_0^{\mathsf T}G_0(w\tens w)_0=\boldsymbol{\ell}.
	\end{equation}
	The entries of $R_0$ and $G_0$ are regular at $q=1$.
	The diagonal entries of $G_0(1)$ are $1,1,9/4,7,8,8$,
	and $R_0(1)$ has rank two by its first and third rows.
	Thus $H(1)$ is positive definite and $\det H\ne0$.
	Since $\Pi_2$ is an orthogonal projection,
	Proposition~\ref{prop:smallspaces}, together with
	\eqref{eq:smallGram} and~\eqref{eq:zeropairings}, gives
	\begin{equation}\label{eq:zeroprojection}
		\Pi_2((w\tens w)_0)
		=R_0(R_0^{\mathsf T}G_0R_0)^{-1}
		R_0^{\mathsf T}G_0(w\tens w)_0
		=R_0H^{-1}\boldsymbol{\ell}.
	\end{equation}
	As $S^2=z^2+2+z^{-2}$, \eqref{eq:restriction},
	\eqref{eq:secondsphere}, \eqref{eq:zeropairings}--\eqref{eq:zeroprojection} yield
	\begin{equation}\label{eq:vGram}
		v+2=q^{22}(f\tens f)\Pi_2((w\tens w)_0)
		=q^{22}\boldsymbol{\ell}^{\mathsf T}H^{-1}\boldsymbol{\ell}.
	\end{equation}
	Applying the two-by-two inverse formula to $H$ in
	\eqref{eq:smallGram}, using \eqref{eq:Rzero} and~\eqref{eq:Gzero},
	reduces \eqref{eq:vGram} to
	\begin{equation}\label{eq:vcleared}
		D(q^{26}-1)v
		=(q^4-1)\bigl[q^{16}A-(q^{10}-1)(q^{24}-1)\bigr].
	\end{equation}
	Substituting
	\[
	\begin{gathered}
		D=q^6[7],\qquad
		q^{26}-1=q^{13}(q-q^{-1})[13],\qquad
		q^4-1=q^2(q-q^{-1})[2],\\
		q^{16}A-(q^{10}-1)(q^{24}-1)
		=q^{17}\bigl([2]-(q-q^{-1})^2[5][12]\bigr)
	\end{gathered}
	\]
	into \eqref{eq:vcleared} gives the formula for $v$ in
	\eqref{eq:uclosed}.
\end{proof}

\section{Reconstruction of the radial operator}
\label{sec:reconstruction}

For $\nu\in\mathsf P\otimes\mathbb Q$, let $d_4(\nu)$ be its
$\alpha_4$-coefficient. Give $\widetilde U$ the
$\mathbb Q/2\mathbb Z$-grading
\[
\deg E_i=\delta_{i4},\quad
\deg F_i=0,\quad
\deg K_\nu=d_4(\nu)\!\!\!\!\!\pmod{2\mathbb Z}.
\]
All defining relations are homogeneous; for the commutators,
$\delta_{i4}\equiv-\delta_{i4}\!\pmod{2\mathbb Z}$.
Write $u_{\overline0}$ for the degree-zero component of $u$.

The black generators have degree zero. Moreover,
$\Theta\nu=\nu$ implies $d_4(\nu)=0$, while $d_4(\gamma)=1$
gives $\deg(E_\gamma K_4^{-1})=0$.
Thus $\widetilde B$ has degree zero, as does $\mathcal N^-$,
generated by the right adjoint black orbit of $F_4$.
Consequently, \eqref{eq:HCdefinition} makes $\HC$ graded.
Let $\widetilde B_+=\ker(\eps|_{\widetilde B})$.
Then the subspace
$J=\widetilde B_+\widetilde U+\widetilde U\widetilde B_+$
is graded as well.

Let $\Rad'$ be the unshifted radial part of
\cite[Lemma~3.1 and Definition~3.2]{Sch26}, and write
$\operatorname{supp}_T$ for translation support.

\begin{lemma}\label{lem:reconstruct}
	Let $C_0\subset\mathbb Z$ be finite.
	An operator
	$R=\sum_{s\in C_0}a_sT_{s\beta}$, with $a_s\in\mathcal F$,
	is uniquely determined by $RS^j$ $(0\leq j<|C_0|)$.
\end{lemma}

\begin{proof}
	If $C_0=\varnothing$, then $R=0$.
	Otherwise, the given values yield the linear system
	\[
	\sum_{s\in C_0}a_s(q^{2s}z+q^{-2s}z^{-1})^j
	=RS^j
	\ (0\leq j<|C_0|).
	\]
	Since $q$ is indeterminate, the Laurent polynomials
	$q^{2s}z+q^{-2s}z^{-1}$ $(s\in C_0)$ have pairwise distinct
	coefficients of $z$, and hence are distinct.
	The coefficient matrix is therefore Vandermonde with nonzero
	determinant in $\kk(z)\subset\mathcal F$.
	Thus the system uniquely determines all $a_s$, and hence $R$.
\end{proof}

Put $\lambda_m=q^{2m}+q^{-22-2m}$ $(m\in\mathbb Z_{\geq0})$.
Extend $\Theta$ to $\mathcal F$ by
$\Theta(e^\mu)=e^{\Theta\mu}$, fixing scalars.
From $\Theta\beta=-\beta$ and $\Theta|_{\beta^\perp}=\id$, we have $\Theta z=z^{-1}$, $\Theta T_{s\beta}\Theta^{-1}=T_{-s\beta}$ $(s\in\mathbb Q)$.
The shifted torus action is
\[
\Theta\mathbin{\cdot}K_\nu
=q^{(\rho,\Theta\nu-\nu)}K_{\Theta\nu}.
\]
In particular,
$\Theta\mathbin{\cdot}K_\beta=q^{-22}K_{-\beta}$ and
$\Theta\mathbin{\cdot}K_{-\beta}=q^{22}K_\beta$,
so $K_\beta+q^{-22}K_{-\beta}$ is invariant under the shifted
action of $\Theta$.

\begin{proposition}\label{prop:radialinput}
	There exists $x\in Z_{\widetilde U}(\widetilde B)$ such that
	\begin{gather*}
		\HC(x)=K_\beta+q^{-22}K_{-\beta},\ 
		\operatorname{supp}_T\Rad(x)\subseteq\{-\beta,0,\beta\},\ 
		\Rad(x)\mathsf{F}_m=\lambda_m\mathsf{F}_m
		\ (m\in\mathbb Z_{\geq0}).
	\end{gather*}
\end{proposition}

\begin{proof}
	By \eqref{eq:sphericallattice} and \cite[Lemma~3.5]{Sch26},
	there is $y\in Z_{\widetilde U}(\widetilde B)$ with
	$\HC(y)=K_\beta+q^{-22}K_{-\beta}$.
Since $\widetilde B$ is concentrated in degree zero, $\HC$ is
graded, and $d_4(\beta)=2$, the element $x=y_{\overline0}$
lies in $Z_{\widetilde U}(\widetilde B)$ and satisfies
$\HC(x)=K_\beta+q^{-22}K_{-\beta}$.
	
	Write $\Rad'(x)=\sum_{\nu\in\mathbb Q\beta}a_\nu T_\nu$
	with finite support.
	For a Zariski-dense set of $h\in\mathsf Y\otimes\mathbb Q$,
	\cite[Lemma~3.1]{Sch26} gives $K_hx-\sum_\nu a_\nu(h)K_hK_\nu\in J$.
	Since $J$ is graded and $a_\nu(h)$ are scalars, taking
	the component of degree $\deg K_h$ gives $K_hx-
	\sum_{d_4(\nu)\in2\mathbb Z}a_\nu(h)K_hK_\nu\in J$.
	
	The truncated sum therefore satisfies the same radial congruence.
	By the uniqueness in \cite[Lemma~3.1]{Sch26}, deleting the
	other summands does not change $\Rad'(x)$.
	Since $d_4(s\beta)=2s$,
	$\operatorname{supp}_T\Rad'(x)\subseteq\mathbb Z\beta$.
	The $-\rho$ shift preserves translation support, so
	$\operatorname{supp}_T\Rad(x)\subseteq\mathbb Z\beta$.
	
	Moreover, \eqref{eq:radialaction} and
	\cite[Lemma~3.5]{Sch26} give
	\begin{equation}\label{eq:radialspectrum}
		\Rad(x)\mathsf{F}_m
		=\HC(x)(m\beta)\mathsf{F}_m
		=\lambda_m\mathsf{F}_m
		\ (m\in\mathbb Z_{\geq0}).
	\end{equation}
	Since $\Theta$ fixes scalars and
	$\Theta\mathsf{F}_m=\mathsf{F}_m$ by
	Lemma~\ref{lem:shiftedsymmetry}, \eqref{eq:radialspectrum} gives
	\[
	\bigl(\Rad(x)-\Theta\Rad(x)\Theta^{-1}\bigr)
	\mathsf{F}_m=0
	\ (m\geq0).
	\]
	The $\mathsf{F}_m$ form a basis of $\kk[S]$, so this difference
	annihilates every $S^j$, $j\geq0$.
	It has finite integral translation support; hence
	Lemma~\ref{lem:reconstruct} yields $\Rad(x)=\Theta\Rad(x)\Theta^{-1}$.

	The upper support bound in the proof of
	\cite[Theorem~4.23]{Sch26}, with leading term $K_\beta$
	in $\HC(x)$, gives 
	\[
	\operatorname{supp}_T\Rad(x)
	\subseteq\beta-\tfrac12\mathbb Z_{\geq0}\beta.
	\]
	Combining this bound with reflection symmetry and integral
	support yields
	\[
	\operatorname{supp}_T\Rad(x)
	\subseteq
	\mathbb Z\beta
	\cap\left(\beta-\tfrac12\mathbb Z_{\geq0}\beta\right)
	\cap\left(-\beta+\tfrac12\mathbb Z_{\geq0}\beta\right)
	=\{-\beta,0,\beta\}.\qedhere
	\]
\end{proof}

\begin{proof}[Proof of Theorem~\ref{thm:main}]
	Let $\mathsf e_j$ $(1\leq j\leq4)$ be the elementary symmetric
	functions in $a,b,c,d$. Equations~\eqref{eq:DAW}--\eqref{eq:AWeigenvalues} give
	\begin{align*}
		D_{\rm AW}S
		&=\delta_1S+(1-q^{-2})(\mathsf e_1-\mathsf e_3),\\
		D_{\rm AW}S^2
		&=\delta_2S^2
		+(1-q^{-4})(\mathsf e_1-q^2\mathsf e_3)S
		+(q^2-1)(1-q^{-4})(\mathsf e_2-\mathsf e_4-1).
	\end{align*}
	
	For the parameters in \eqref{eq:AWparameters},
	the formula for $v$ in \eqref{eq:uclosed} gives
	\[
		D(q^{26}-1)v
		=(q^2-1)\bigl[D(\mathsf e_2-q^{24}-1)
		+q^4A^4B_0^2(1+q^{12})\bigr].
	\]
	Together with \eqref{eq:uclosed}, this identity yields
	\begin{equation}\label{eq:uvAW}
\begin{aligned}
 u=\frac{(q^2+1)(\mathsf e_1-q^2\mathsf e_3)}
          {\mathsf e_4q^4-1},\ 
 (\mathsf e_4q^2-1)v
 =(q^2-1)(\mathsf e_2-\mathsf e_4-1)
   +\frac{uq^2}{q^2+1}(\mathsf e_1-\mathsf e_3).
\end{aligned}
\end{equation}
	Proposition~\ref{prop:second}, together with \eqref{eq:firstAW},
\eqref{eq:firstsphere}, and~\eqref{eq:uvAW}, gives
$D_{\rm AW}\mathsf{F}_j=\delta_j\mathsf{F}_j$ $(j=0,1,2)$.

	Choose $x$ as in Proposition~\ref{prop:radialinput}.
	Since $\lambda_m=q^{-22}\delta_m+1+q^{-22}$,
	the operators $\Rad(x)$ and
	$q^{-22}D_{\rm AW}+(1+q^{-22})\id$ agree on $\mathsf{F}_0,\mathsf{F}_1,\mathsf{F}_2$.
	Both have support in $\{-\beta,0,\beta\}$, and these polynomials
	span the polynomials of degree at most two in $S$.
	Lemma~\ref{lem:reconstruct}, with $C_0=\{-1,0,1\}$,
	therefore proves \eqref{eq:mainradial}.
	
	Equation~\eqref{eq:radialspectrum} now gives
	$D_{\rm AW}\mathsf{F}_m=\delta_m\mathsf{F}_m$ for every $m\geq0$.
	Lemma~\ref{lem:shiftedsymmetry} and the defining uniqueness
	of $P_m$ imply $\mathsf{F}_m=P_m$, proving \eqref{eq:main}
and hence \eqref{eq:Schconjecture}.
\end{proof}

\section{Parameter transport and normalizations}
\label{sec:transport}

The general torus-equivariance underlying parameter transport is
established in \cite[Theorem~6.3]{Let03}.
We give a direct formula in the present conventions.
Fix $\delta,\sigma\in\kk^\times$ and put $\eta=\delta/\sigma$.

\begin{corollary}\label{cor:transport}
For every $m\in\mathbb Z_{\geq0}$,
\begin{equation}\label{eq:transport}
 \Res\phi_m^{\delta,\sigma}(z)
 =q^{-11m}\eta^mP_m(z/\eta).
\end{equation}
If $\delta,\sigma\in\{q^5,-q^5\}$, then
\begin{equation}\label{eq:signedAW}
 \Res\phi_m^{\delta,\sigma}(z)
 =q^{-11m}P_m(z;\eta a,\eta b,\eta c,\eta d\mid q^2).
\end{equation}
\end{corollary}
\begin{proof}
	The fixed-torus conditions force every spherical vector and covector
	to be supported on
	$(\mathsf Y^\Theta)^\perp\cap\mathsf P=\mathbb Z\beta$.
	Define
	\[
	w_m^\sigma
	=\sum_{k\in\mathbb Z}
	\left(\frac{\sigma}{q^5}\right)^{k-m}(w_m)_{k\beta},
	\quad
	f_m^\delta\big|_{V(m\beta)_{k\beta}}
	=\left(\frac{\delta}{q^5}\right)^{m-k}
	f_m\big|_{V(m\beta)_{k\beta}},
	\]
	with $f_m^\delta=0$ on the other weight spaces.
	These rescalings preserve invariance under $U_X$ and
	$U^{0,\Theta}$.
	Lemma~\ref{lem:allspheres} gives $B_4w_m=0$ and $f_mB_4=0$.
	For $u\in V(m\beta)_{k\beta+\alpha_4}$, the definitions of
	$w_m^\sigma$ and $f_m^\delta$, together with
	$\gamma+\alpha_4=\beta$, yield
	\begin{align*}
		\bigl((F_4+\sigma E_\gamma K_4^{-1})w_m^\sigma\bigr)_{k\beta-\alpha_4}
		&=\left(\frac{\sigma}{q^5}\right)^{k-m}
		(B_4w_m)_{k\beta-\alpha_4}=0,\\
		f_m^\delta\bigl((F_4+\delta E_\gamma K_4^{-1})u\bigr)
		&=\left(\frac{\delta}{q^5}\right)^{m-k}f_m(B_4u)=0.
	\end{align*}
	Other input weights give zero by support.
	The inverse rescalings identify the invariant spaces with the
	standard spherical lines and preserve highest-weight normalization,
	proving uniqueness.
	
	The factors at weight $k\beta$ have product $\eta^{m-k}$, so
	\eqref{eq:restriction} gives
	\begin{equation}\label{eq:restrictiontransport}
	\Res\phi_m^{\delta,\sigma}(z)
	=\eta^m(\Res\phi_m)(z/\eta).
	\end{equation}
Substituting \eqref{eq:main} into \eqref{eq:restrictiontransport} gives
\eqref{eq:transport}. 
	
	For $\eta^2=1$, replacing $z$ by $z/\eta$ in
	\eqref{eq:AWcoefficient} changes the parameters to
	$(\eta a,\eta b,\eta c,\eta d)$, preserving $abcd$ and the
	eigenvalues in \eqref{eq:AWeigenvalues}.
	Evaluating the eigenvalue equation for $P_m$ at $z/\eta$
	shows that $\eta^mP_m(z/\eta)$ is an eigenpolynomial for these
	parameters with eigenvalue $\delta_m$.
	Its leading term in $S$ is $\eta^{2m}S^m=S^m$.
	Uniqueness and \eqref{eq:transport} give \eqref{eq:signedAW}.
\end{proof}

For $m=1$ and $(\delta,\sigma)=(-q^5,q^5)$,
\eqref{eq:transport} gives
\[
 q^{11}\Res\phi_1^{-q^5,q^5}
 =-P_1(-z)=S-\frac{[2]}{[3]},
 \quad P_1(z)=S+\frac{[2]}{[3]}.
\]
Since $[2]/[3]\ne0$, these monic polynomials are not proportional,
illustrating the dependence on the relative sign.

The positivity argument in the proof of
Lemma~\ref{lem:shiftedsymmetry} gives $(w_m,w_m)\ne0$
for every $m\ge0$.
Equations~\eqref{eq:restriction} and~\eqref{eq:main} then yield
\[
\phi_m(1)
=(w_m,w_m)=q^{-11m}P_m(q^{11})\ne0,\quad
\Res\!\left(\frac{\phi_m}{\phi_m(1)}\right)
=\frac{P_m(z)}{P_m(q^{11})}.
\]
\section*{Acknowledgements}
The author thanks Philip Schl\"osser for helpful comments on the manuscript
and for clarifying the parameter convention in his FII conjecture.


\begin{thebibliography}{99}
\raggedright
\bibitem{Kol14} S. Kolb,
\emph{Quantum symmetric Kac--Moody pairs},
Adv. Math. \textbf{267} (2014), 395--469.

\bibitem{Let03} G. Letzter,
\emph{Quantum symmetric pairs and their zonal spherical functions},
Transform. Groups \textbf{8} (2003), 261--292.

\bibitem{Let04} G. Letzter,
\emph{Quantum zonal spherical functions and Macdonald polynomials},
Adv. Math. \textbf{189} (2004), 88--147.

\bibitem{Lus} G. Lusztig,
\emph{Introduction to Quantum Groups},
Modern Birkh\"auser Classics,
Birkh\"auser Boston, Boston, MA, 2010.

\bibitem{Mee26} S. Meereboer,
\emph{Symmetries for spherical functions of type $\chi$ for quantum
symmetric pairs}, Transform. Groups \textbf{31} (2026), 1851--1894.

\bibitem{MS26} S. Meereboer and P. Schl\"osser,
\emph{Quantum matrix-spherical functions as orthogonal polynomials},
arXiv:2511.23367v2 (2026).

\bibitem{Sch26} P. Schl\"osser,
\emph{Zonal Spherical Functions of Quantum Symmetric Pairs},
arXiv:2606.27901v1 (2026).

\bibitem{Wat23} H. Watanabe,
\emph{Stability of $\imath$canonical bases of irreducible finite type of
real rank one}, Represent. Theory \textbf{27} (2023), 1--29.

\end{thebibliography}
\end{document}